\documentclass[12pt]{amsart}
\usepackage{amssymb}
\usepackage{latexsym}
\usepackage{pb-diagram}
\usepackage{lamsarrow}
\usepackage{pb-lams}

\newtheorem{theorem}{Theorem}[section]

\newtheorem{definition}[theorem]{Definition}

\newtheorem{lemma}[theorem]{Lemma}

\begin{document}

\title[The Jensen covering property]{The Jensen covering property}
\author{E. Schimmerling}
\author{W.H. Woodin}
\date{March 13, 1996}
\address{Department of Mathematics, University of California, 
Irvine, CA 92697-4065, USA}
\email{eschimme@math.uci.edu}
\address{Department of Mathematics, University of California,
Berkeley, California 94720, USA}
\email{woodin@math.berkeley.edu}
\thanks{The research of the first author 
was partially supported by an NSF 
Mathematical Sciences Postdoctoral Research Fellowship at the 
Massachusetts Institute of Technology.}
\thanks{The research of the second author was partially supported by the
NSF.}
\maketitle

\newcommand\res{{{\upharpoonright}}}
\newcommand\exist{\exists}
\newcommand\is{\unlhd}
\newcommand\pis{\lhd}
\newcommand\npis{\ntriangleleft}


\newcommand\OR{{\mathrm{OR}}}
\newcommand\ult{{\mathrm{ult}}}
\newcommand\pred{{\mathrm{pred} \,}}
\newcommand\rt{{\mathrm{root} \,}}
\newcommand\decap{{\mathrm{decap}}}
\newcommand\otp{{\mathrm{otp}}}


\newcommand\avec{{\vec a}}
\newcommand\bvec{{\vec b}}
\newcommand\cvec{{\vec c}}
\newcommand\dvec{{\vec d}}
\newcommand\kvec{{\vec k}}
\newcommand\mvec{{\vec m}}
\newcommand\nvec{{\vec n}}
\newcommand\uvec{{\vec u}}
\newcommand\vvec{{\vec v}}
\newcommand\xvec{{\vec x}}
\newcommand\yvec{{\vec y}}

\newcommand\Avec{{\vec A}}
\newcommand\Bvec{{\vec B}}
\newcommand\Cvec{{\vec C}}
\newcommand\Dvec{{\vec D}}
\newcommand\Evec{{\vec E}}

\newcommand\gavec{{{\vec {\ga} }}}
\newcommand\gbvec{{{\vec {\gb} }}}
\newcommand\gkvec{{{\vec {\gk} }}}
\newcommand\etavec{{{\vec {\eta} }}}
\newcommand\glvec{{{\vec {\gl} }}}
\newcommand\gmvec{{{\vec {\mu} }}}
\newcommand\gnvec{{{\vec {\nu} }}}
\newcommand\gsvec{{{\vec {\gs} }}}

\newcommand\gLvec{{{\vec {\gL} }}}
\newcommand\gOvec{{{\vec {\gO} }}}

\newcommand\cMvec{{{\vec {\mathcal M} }}}
\newcommand\cNvec{{{\vec {\mathcal N} }}}
\newcommand\cPvec{{{\vec {\mathcal P} }}}
\newcommand\cQvec{{{\vec {\mathcal Q} }}}
\newcommand\cRvec{{{\vec {\mathcal R} }}}
\newcommand\cSvec{{{\vec {\mathcal S} }}}


\newcommand\htil{{\widetilde h}}
\newcommand\ktil{{\widetilde k}}

\newcommand\Ftil{{\widetilde F}}
\newcommand\Gtil{{\widetilde G}}
\newcommand\Ktil{{\widetilde K}}

\newcommand\cPtil{{\widetilde {\cP}}}
\newcommand\cQtil{{\widetilde {\cQ}}}
\newcommand\cTtil{{\widetilde {\cT}}}
\newcommand\cUtil{{\widetilde {\cU}}}
\newcommand\cVtil{{\widetilde {\cV}}}

\newcommand\gbtil{{\widetilde {\gb}}}
\newcommand\gktil{{\widetilde {\gk}}}
\newcommand\gptil{{\widetilde {\gp}}}
\newcommand\gdtil{{\widetilde {\gd}}}
\newcommand\gstil{{\widetilde {\gs}}}
\newcommand\gftil{{\widetilde {\gf}}}
\newcommand\etatil{{\widetilde {\eta}}}
\newcommand\psitil{{\widetilde {\psi}}}

\newcommand\gOtil{{\widetilde {\gO}}}


\newcommand\sdot{{\dot s}}
\newcommand\xdot{{\dot x}}
\newcommand\Edot{{\dot E}}
\newcommand\Fdot{{\dot F}}
\newcommand\ggdot{{\dot \gg}}
\newcommand\gmdot{{\dot \gm}}
\newcommand\gndot{{\dot \gn}}
\newcommand\gtdot{{\dot \gt}}


\newcommand\abar{{\overline a}}
\newcommand\bbar{{\overline b}}
\newcommand\cbar{{\overline c}}
\newcommand\ibar{{\overline i}}
\newcommand\pbar{{\overline p}}
\newcommand\qbar{{\overline q}}
\newcommand\rbar{{\overline r}}
\newcommand\sbar{{\overline s}}
\newcommand\tbar{{\overline t}}
\newcommand\ubar{{\overline u}}
\newcommand\vbar{{\overline v}}
\newcommand\wbar{{\overline w}}
\newcommand\xbar{{\overline x}}
\newcommand\ybar{{\overline y}}

\newcommand\Abar{{\overline A}}
\newcommand\Bbar{{\overline B}}
\newcommand\Cbar{{\overline C}}
\newcommand\Dbar{{\overline D}}
\newcommand\Ebar{{\overline E}}
\newcommand\Fbar{{\overline F}}
\newcommand\Kbar{{\overline K}}
\newcommand\Tbar{{\overline T}}
\newcommand\Wbar{{\overline W}}

\newcommand\cHbar{{\overline {\cH}}}
\newcommand\cIbar{{\overline {\cI}}}
\newcommand\cMbar{{\overline {\cM}}}
\newcommand\cNbar{{\overline {\cN}}}
\newcommand\cPbar{{\overline {\cP}}}
\newcommand\cQbar{{\overline {\cQ}}}
\newcommand\cRbar{{\overline {\cR}}}
\newcommand\cSbar{{\overline {\cS}}}
\newcommand\cTbar{{\overline {\cT}}}
\newcommand\cUbar{{\overline {\cU}}}
\newcommand\cVbar{{\overline {\cV}}}

\newcommand\gabar{{\overline \ga}}
\newcommand\gbbar{{\overline \gb}}
\newcommand\ggbar{{\overline \gg}}
\newcommand\gdbar{{\overline \gd}}
\newcommand\gkbar{{\overline \gk}}
\newcommand\glbar{{\overline \gl}}
\newcommand\gnbar{{\overline \gn}}
\newcommand\gpbar{{\overline \gp}}
\newcommand\grbar{{\overline \gr}}
\newcommand\gsbar{{\overline \gs}}
\newcommand\gtbar{{\overline \gt}}

\newcommand\gGbar{{\overline \gG}}
\newcommand\gObar{{\overline \gO}}
\newcommand\gUbar{{\overline \gU}}


\newcommand\fA{\mathfrak A}
\newcommand\fB{\mathfrak B}
\newcommand\fC{\mathfrak C}
\newcommand\fM{\mathfrak M}
\newcommand\fR{\mathfrak R}
\newcommand\fS{\mathfrak S}
\newcommand\fT{\mathfrak T}


\newcommand{\cA}{{\mathcal A}}
\newcommand{\cB}{{\mathcal B}}
\newcommand{\cC}{{\mathcal C}}
\newcommand{\cD}{{\mathcal D}}
\newcommand{\cE}{{\mathcal E}}
\newcommand{\cF}{{\mathcal F}}
\newcommand{\cG}{{\mathcal G}}
\newcommand{\cH}{{\mathcal H}}
\newcommand{\cI}{{\mathcal I}}
\newcommand{\cJ}{{\mathcal J}}
\newcommand{\cK}{{\mathcal K}}
\newcommand{\cL}{{\mathcal L}}
\newcommand{\cM}{{\mathcal M}}
\newcommand{\cN}{{\mathcal N}}
\newcommand{\cO}{{\mathcal O}}
\newcommand{\cP}{{\mathcal P}}
\newcommand{\cQ}{{\mathcal Q}}
\newcommand{\cR}{{\mathcal R}}
\newcommand{\cS}{{\mathcal S}}
\newcommand{\cT}{{\mathcal T}}
\newcommand{\cU}{{\mathcal U}}
\newcommand{\cV}{{\mathcal V}}
\newcommand{\cW}{{\mathcal W}}
\newcommand{\cX}{{\mathcal X}}
\newcommand{\cY}{{\mathcal Y}}
\newcommand{\cZ}{{\mathcal Z}}


\newcommand\lb{\lbrack}
\newcommand\rb{\rbrack}
\newcommand\lcb{\lbrace}
\newcommand\rcb{\rbrace}


\newcommand\ra{\rightarrow}
\newcommand\la{\leftarrow}

\newcommand\thra{\twoheadrightarrow}

\newcommand\lra{\longrightarrow}
\newcommand\lla{\longleftarrow}

\newcommand\llra{\longleftrightarrow}

\newcommand\Ra{\Rightarrow}
\newcommand\La{\Leftarrow}

\newcommand\Lra{\Longrightarrow}
\newcommand\Lla{\Longleftarrow}

\newcommand\Llra{\Longleftrightarrow}
\renewcommand\iff{\Leftrightarrow}

\newcommand{\HOD}{\mathrm{HOD}}
\newcommand{\AD}{\mathrm{AD}}

%
%

\newcommand{\ga}{\alpha}     
\newcommand{\gb}{\beta}      
\renewcommand{\gg}{\gamma}   
\newcommand{\gd}{\delta}     
\renewcommand\ge{\varepsilon}
\newcommand{\gz}{\zeta}      
\newcommand{\gth}{\theta}    
\newcommand{\gi}{\iota}      
\newcommand{\gk}{\kappa}  
\newcommand{\gl}{\lambda}    
\newcommand{\gm}{\mu}        
\newcommand{\gn}{\nu}        
\newcommand{\gx}{\xi}        
\newcommand{\gom}{\omicron}  
\newcommand{\gp}{\pi}        
\newcommand{\gr}{\rho}       
\newcommand{\gs}{\sigma}     
\newcommand{\gt}{\tau}       
\newcommand{\gu}{\upsilon}   
\newcommand{\gph}{\varphi}      
\newcommand{\gf}{\varphi}      
\newcommand{\gch}{\chi}      
\newcommand{\gps}{\psi}      
\newcommand{\go}{\omega}     

\newcommand{\gA}{A}     
\newcommand{\gB}{B}      
\newcommand{\gG}{\Gamma}     
\newcommand{\gD}{\Delta}     
\newcommand{\gEp}{E}  
\newcommand{\gZ}{Z}      
\newcommand{\gEe}{H}      
\newcommand{\gTh}{\Theta}    
\newcommand{\gI}{I}      
\newcommand{\gK}{K}     
\newcommand{\gL}{\Lambda}    
\newcommand{\gM}{M}        
\newcommand{\gN}{N}        
\newcommand{\gX}{\Xi}        
\newcommand{\gOm}{O}  
\newcommand{\gP}{\Pi}        
\newcommand{\gR}{P}       
\newcommand{\gS}{\Sigma}     
\newcommand{\gT}{T}       
\newcommand{\gU}{\Upsilon}   
\newcommand{\gPh}{\Phi}      
\newcommand{\gCh}{X}      
\newcommand{\gPs}{\Psi}      
\newcommand{\gO}{\Omega}         

\newcommand{\bA}{{\bf A}}  
\newcommand{\bB}{{\bf B}}      
\newcommand{\bG}{\boldGamma}     
\newcommand{\bD}{\boldDelta}     
\newcommand{\bEp}{{\bf E}}  
\newcommand{\bZ}{{\bf Z}}      
\newcommand{\bEe}{{\bf H}}      
\newcommand{\bTh}{\boldTheta}    
\newcommand{\bI}{{\bf I}}     
\newcommand{\bK}{{\bf K}}     
\newcommand{\bL}{{\bf L}}    
\newcommand{\bM}{{\bf M}}        
\newcommand{\bN}{{\bf N}}        
\newcommand{\bX}{\boldXi}        
\newcommand{\bOm}{{\bf O}}  
\newcommand{\bP}{\boldPi}        
\newcommand{\bR}{{\bf P}}       
\newcommand{\bS}{\boldSigma}     
\newcommand{\bT}{{\bf T}}       
\newcommand{\bU}{\boldUpsilon}   
\newcommand{\bPh}{\boldPhi}      
\newcommand{\bCh}{{\bf X}}      
\newcommand{\bPs}{\boldPsi}      
\newcommand{\bO}{\boldOmega}     

%
%

\newcommand{\ha}{\aleph}
\newcommand{\hb}{\beth}
\newcommand{\hg}{\gimel}
\newcommand{\hd}{\daleth}

%
%

\newcommand{\forces}{\Vdash}
\newcommand{\decides}{\parallel}

\newcommand{\dom}{\mathrm{dom}}
\newcommand{\ran}{\mathrm{ran}}
\newcommand{\crit}{\mathrm{crit}}
\newcommand{\card}{\mathrm{card}}
\newcommand{\supp}{\mathrm{supp}}
\newcommand{\support}{\mathrm{support}}
\newcommand{\cf}{\mathrm{cf}}
\newcommand{\lh}{\mathrm{lh}}
\newcommand{\rank}{\mathrm{rank}}
\newcommand{\id}{\mathrm{id}}
\newcommand{\ot}{\mathrm{o.t.}}
\section{Introduction}

The Jensen covering lemma says that either $L$ has a club class of
indiscernibles, or else,
for every uncountable set $A$ of ordinals, 
there is a set $B \in L$ with $A \subseteq B$
and $\card (B) = \card (A)$.  
One might hope to extend Jensen's covering lemma to richer weasels, 
that is, to inner models of the form
$L[\Evec ]$ where $\Evec$ is a coherent sequence of extenders of the kind
studied in Mitchell-Steel \cite{MiSt}.
The papers \cite{MiSt}, \cite{many}, and \cite{tame} 
show how to construct weasels with Woodin cardinals and more.
But, as Prikry forcing shows,
one cannot expect too direct a generalization
of Jensen's covering lemma to weasels with measurable cardinals.

Recall from \cite{MiSt} that if 
$L[\Evec ]$ is a weasel and $\ga$ is an ordinal, 
then either $E_\ga = \emptyset$, or else
$E_\ga$ is an extender over $J^\Evec_\ga$.
Let us say that $L[\Evec ]$ is a {\bf lower-part} weasel iff
for every ordinal $\ga$, 
$E_\ga$ is not a total extender over $L[\Evec ]$.
In other words, if $L[\Evec ]$ is a lower-part weasel,
then no cardinal in $L[\Evec ]$ is measurable as witnessed by an
extender on $\Evec$.
But a lower-part weasel $L[\Evec ]$ could be rich in the sense that
it may have levels $J^\Evec_\ga$ satisfying 
ZFC + ``there are many Woodin cardinal''.
We do not impose any bounds on the large cardinal axioms true
in the levels of a lower-part weasel.

In this paper we show that if $L[\Evec ]$ is an iterable,
lower-part weasel,
then either $L[\Evec ]$ has indiscernibles,
or else $L[\Evec ]$ satisfies the Jensen covering property.
The iterability and indiscernibility 
that we mean will be made precise in due course. 
Our result says, in a new way, that
Prikry forcing is the essential limitation on
extensions of Jensen's covering lemma.

This paper is organized as follows.  In \S\S2-5, we prove several
theorems regarding successor cardinals in lower-part weasels.  
Those in \S2 and \S3 might be seen as corollaries to the proofs of 
the main theorems in \cite{MiSchSt} and \cite{MiSch} respectively.
Iterability is discussed in \S4, where we sketch why
internal strategic iterability for countable
realizable premice suffices for all of the results in this paper.
In \S5, we prove our extension to Jensen's covering lemma.

That Theorem~2.4 is a corollary to \cite{MiSchSt} was observed by Woodin
\cite{W}, who used it to define a core model based on the failure of
$L(\mathbb R)$-determinacy.
We were originally and primarily motivated by this 
application and give a rough description of it at the end of \S2.

The results presented in this paper build directly on the earlier work
of Mitchell, Schimmerling, and Steel.

\section{Successor cardinals with countable closure}

We assume that the reader is familiar with \cite{MiSt} and \cite{cmip},
and adopt the notation established there.
Throughout the paper,
if $\gp$ is an embedding, then we shall write $\cE (\gp )$, instead of
the more common notation $E_\gp$,
for the long extender derived from $\gp$.
For the purposes of this paper, an iteration tree is a normal, padded
iteration tree.
If $\cT$ is an iteration tree of successor length, then we shall write
$\cM_\infty ( \cT )$ for the final model of $\cT$.
If $\cM$ and $\cN$ are premice,
then $\cM \is \cN$ means that $\cM$ is an initial segment of
$\cN$, and
$\cM \pis \cN$ means that $\cM$ is a proper initial segment of $\cN$.

\begin{definition}\label{2.1}
A premouse $\cM$ is a {\bf lower-part} premouse if
and only if $\cM$ is passive and
$\Edot^\cM_\ga$ is not a total extender over $\cM$
whenever $\ga < \OR^\cM$.
\end{definition}

\begin{definition}\label{2.2}
Let $\cM$ be a premouse.
Then $\cM$ is {\bf iterable} iff
$\cM$ is $\ga$-iterable for every ordinal $\ga$ in the sense of
\cite[Definition~2.8]{cmip}.
That is, $\cM$ is $\ga$-iterable iff player II has a winning strategy in the
iteration game $\cG^* (\cM , \ga )$ of length $\ga$ on $\cM$.
If $\ga$ is an ordinal and $\gk < \OR^\cM$,
then $\cM$ is {\bf $\ga$-iterable above $\gk$} iff
player II has a winning strategy in 
the iteration game of length $\ga$ on $\cM$ where 
player I is restricted to playing extenders with critical points $\geq \gk$.
$\cM$ is {\bf iterable above $\gk$} iff
$\cM$ is $\ga$-iterable above $\gk$ for all ordinals $\ga$.
\end{definition}

\begin{definition}\label{2.3}
$\cM$ is an {\bf $x$-premouse}
iff $x$ is a set of ordinals and 
$\cM$ is a premouse built over $x$.  
\end{definition}

Definition~2.3 means, among other things,
that if $E$ is an extender from the $\cM$-sequence,
then $\cup x < \crit (E)$.
All definitions involving premice have the obvious extensions to
$x$-premice.

Recall that a cardinal $\gk$ is countably closed iff
$\gm^{\ha_0} < \gk$ whenever $\gm < \gk$.

\begin{theorem}\label{2.4}
Suppose that $W = L[\Evec , x]$ is an iterable, lower-part $x$-weasel.
Let $\gk$ be a cardinal of $W$ such that $\cup x < \card ( \gk )$ and 
$\card (\gk )$ is a countably closed cardinal. 
Put $\gl = ( \gk^+ )^W$ and suppose that $$\cf (\gl ) < \card (\gk ).$$
Then, given $\gd < \card (\gk )$, there exists $\gm \geq \gd$ and
a normal, countably complete ultrafilter $U$ on $\cP (\gm ) \cap W$ such
that $U$ is weakly amenable to $W$.
\end{theorem}

\begin{definition}\label{2.5}
An iteration tree $\cT$ is {\bf trivial} iff
$E^\cT_\eta = \emptyset$ whenever $\eta < \lh ( \cT)$.  
That is, $\cT$ is trivial iff $\cT$ 
consists exclusively of padding.
\end{definition}

\begin{definition}\label{2.6}
An iteration tree $\cT$ is {\bf dropping} iff
$\eta + 1 \in \cD^\cT$ whenever $0 = \pred^\cT (\eta +1 )$.
\end{definition}

\begin{lemma}\label{2.7}
If $\cT$ is an iteration tree on a lower part-part premouse
and $E^\cT_0 \not= \emptyset$, then $\cT$ is dropping.
\end{lemma}

\begin{proof}[Proof of Lemma 2.7]
Since $\cM^\cT_0$ is a lower-part premouse, we have that $1 \in \cD^\cT$.
Suppose that $\pred^\cT (\eta +1) = 0 < \eta$ but $\eta +1 \not\in \cD^\cT$.
Then a simple calculation shows that
$$\crit ( E^\cT_\eta ) < \gr_\go \left( (\cM_1^* )^\cT \right) 
\leq \crit ( E^\cT_0 ) < \lh ( E^\cT_0 ) < \lh (
E^\cT_\eta ).$$
So $\crit (E^\cT_0 )$ is a limit of generators of $E^\cT_\eta$.
By the initial segment condition on good extender sequences,
there is are initial segments of
$E^\cT_\eta \res \crit (E^\cT_0 )$ on the $\cM^\cT_0$-sequence,
and these initial segments are total extenders over $\cM^\cT_0$.
But $\cM^\cT_0$ is a lower-part premouse, a contradiction.
\end{proof}

\begin{definition}\label{2.8}
An iteration tree $\cT$ is {\bf thorough} iff 
$E^\cT_\eta \not= \emptyset$ and
$$\cM^\cT_\eta \parallel \lh (E^\cT_\eta )$$ is a lower part
premouse whenever $\eta < \lh (\cT )$.
\end{definition}

\begin{definition}\label{2.9}
An extender $E$ is a {\bf measure} iff the only
generator of $E$ is $\crit (E)$.
\end{definition}

\begin{lemma}\label{2.10}
Suppose that $\cT$ is a thorough iteration tree.
Then $E^\cT_\eta$ is a measure whenever $\eta < \lh (\cT )$,
so $\cT$ is linear.
Moreover, $\cT$ has the following very simple form:
$$(\gg + 1 , \eta ]_T \cap \cD^\cT = \emptyset
\ \Lra \
E^\cT_\eta = i^\cT_{\gg , \eta} (E^\cT_\gg )$$
whenever $\gg < \eta < \lh (\cT)$.
\end{lemma}

\begin{proof}[Proof of Lemma 2.10]
The first assertion is clear from the initial segment condition on
$\cM^\cT_\eta$ and the thoroughness of $\cT$.
Suppose that $\gg < \eta < \lh (\cT)$.
If 
$$\lh (E^\cT_\eta ) 
> \lh \left( i^\cT_{\gg , \eta} (E^\cT_\gg ) \right),$$
then $i^\cT_{\gg , \eta} (E^\cT_\gg )$
witnesses that $\cM^\cT_\eta \parallel \lh (E^\cT_\eta )$ is not a
lower-part premouse.
Suppose that $$\lh (E^\cT_\eta ) 
< \lh \left( i^\cT_{\gg , \eta} (E^\cT_\gg ) \right).$$
Then $E^\cT_\eta$ witnesses that
$$i^\cT_{\gg , \eta} \left(\cM^\cT_\gg \parallel \lh (E^\cT_\gg )
\right)$$
is not a lower-part premouse.
By elementarity,
$\cM^\cT_\gg \parallel \lh (E^\cT_\gg )$ is not a lower-part premouse,
a contradiction.
\end{proof}

\begin{proof}[Proof of Theorem~2.4]
For the reader familiar with \cite{MiSchSt}, we remark that
this proof begins just like the proof of the main result there.
But the current proof is substantially simpler, as there
is no need for an inductive proof that the lift-ups are strong.

Let $\gO$ be a regular cardinal $>\gl$, and let $X$ be an elementary
substructure of $V_\gO$ with
$$\card (X) < \card (\gk ),$$
$${}^\go X \subset X,$$
$$\sup (X \cap \gl ) = \gl ,$$
and
$$(\cup x) \cup \gd \cup \{ \cJ^W_\gl \} \subset X.$$
Let $\gp : N \lra V_\gO$ be the inverse of the transitive collapse of $X$.
Say 
$\gp (\gkbar ) = \gk$,
$\gp ( \glbar ) = \gl$,
$\gp ( \Wbar ) = \cJ^W_\gl$,
and 
$\ge = \crit ( \gp )$.  Then $\ge < \gk$.

If $W$ were $1$-small, then we could quote \cite[Theorem~8.2]{MiSt} 
and conclude that $\Wbar$ and $W$ agree below $(\ge^+ )^\Wbar$.
But $W$ is not assumed to satisfy a smallness condition.
Nevertheless, we can argue as follows.
Suppose that $\cM \pis \Wbar$ with $\gr_\go (\cM ) = \ge$.
Let $(\cSbar , \cS)$ be a successful coiteration of $(\cM , \gp ( \cM ))$.
All the models on $\cSbar$ and $\cS$ agree below $\ge$,
and $\ge$ is a cardinal in all the models on $\cSbar$ and $\cS$ except possibly
$\gp ( \cM )$.
Using the initial segment condition on good extender sequences and 
the fact that $\Wbar$ is a lower-part $x$-premouse,
one concludes that all extenders used on the iteration trees $\cSbar$ and
$\cS$ have critical points at least $\ge$.  
Standard arguments using the soundness of $\cM$ and $\gp ( \cM )$ 
now show that both $\cSbar$ and $\cS$ are trivial,
so that $\cM \pis \gp ( \cM )$.
So every such $\cM$ is a proper initial segment of $W$.
Since the initial segments $\cM$ of $\Wbar$ 
with $\gr_\go (\cM ) = \ge$ appear cofinally
in $\Wbar \parallel (\ge^+ )^\Wbar$, 
we have shown that $\Wbar$ and $W$ agree below $(\ge^+ )^\Wbar$.

First, consider the case in which $(\ge^+ )^\Wbar = (\ge^+ )^W$,
so that $\Wbar$ and $W$ have the same subsets of $\ge$.
Then Theorem~2.4 holds with $\gm = \ge$ and
$U$ equal to the measure derived from $\gp$, that is,
$$U = \left\lbrace 
x \in |\Wbar| \mid x \subseteq \ge \ \land \ \ge \in \gp (x)
\right\rbrace.$$
Thus, we may assume without loss of generality that
$(\ge^+ )^\Wbar < (\ge^+)^W$.

Let $(\cTbar , \cT)$ be a successful coiteration of $(\Wbar , W)$,
and say this coiteration has length $\theta +1$.
Note that both $\Wbar$ and $W$ are both lower-part $x$-premice,
and models of ZF$^-$.
Moreover,
$(\ge^+ )^\Wbar \leq \lh (E^\cT_0 ) < (\ge^+ )^W$.

One of $\cTbar$ and $\cT$ must be trivial,
as otherwise 
$$[0,\theta ]_\Tbar \cap \cD^\cTbar \not=\emptyset$$ and 
$$[0,\theta ]_T \cap \cD^\cT \not=\emptyset ,$$
which leads to the usual contradiction as in the comparison lemma
\cite[Lemma~7.1]{MiSt}.
Since $\Wbar$ has fewer subsets of $\ge$ than $W$,
it must be that $\cTbar$ is trivial.
Therefore $\Wbar \pis \cM^\cT_\theta$
and $\cT$ does not involve any padding.
By Lemma~2.7, since $W$ is a lower-part $x$-weasel,
$\cT$ is dropping.
Since $\Wbar$ is a lower-part $x$-premouse,
$\cT$ is thorough.
Thus, $\cT$ has the special form described in Lemma~2.10.

Let $\eta$ be the least ordinal such that $\gn^\cT_\eta > \gkbar$,
if such an ordinal exists, and let $\eta = \theta$ otherwise.
We remark that because $\cT$ is thorough and $\gkbar$ is a cardinal in $\Wbar$,
$$\gn^\cT_\gg > \gkbar \ \Llra \crit(E^\cT_\gg ) \geq \gkbar$$
for all $\gg \leq \theta$.

Let $\cP$ be the longest initial segment of $\cM^\cT_\eta$ with exactly
those subsets of $\gkbar$ that are in $\Wbar$.
Then $\cP$ is a $\gkbar$-sound $x$-premouse.

Let $\cR = \ult (\cP , \cE ( \gp ) \res \gk )$
and $\gptil : \cP \lra \cR$ be the ultrapower map.
Recall that this ultrapower is defined using coordinates 
$a \in [\gk ]^{<\go}$ and certain $\cP$-definable functions 
$f : [\gk ]^{<\go } \lra |\cP |$.
The usual argument, using the countable completeness of $\cE ( \gp ) \res \gk$,
shows that $\cR$ is wellfounded.
By the standard fine structural results in \cite[\S2]{MiSchSt},
$\cR$ and $W$ agree below $\gl$,
$$\gp \res \gl = \gptil \res \gl ,$$
$$(\gk^+ )^\cR = \sup (\gptil'' \glbar ) = \sup (\gp'' \glbar )
= \gp (\glbar ) = \gl,$$
and $\cR$ is $\gk$-sound.

We claim that if $\cP$ is active,
then $\gmdot^\cP \geq \gkbar$.
For, if $\gmdot^\cP < \gkbar$,
then initial segments of $\Fdot^\cP \res \gkbar$ are on the $\Wbar$-sequence,
contradicting that $\Wbar$ is a lower-part $x$-premouse.

The previous claim implies that $\gptil$ is continuous at $(\gmdot^+ )^\cP$,
and hence that $\Fdot^\cR$ measures all sets in $\cR$. 
It follows that $\cR$ is an $x$-premouse (rather than just an
$x$-protomouse, in the terminology of \cite{MiSchSt}).

Countable completeness can also be used, in the standard way,
to show that $\cR$ is iterable above $\gk$.
Note that any coiteration of $(W , \cR )$ uses extenders with critical points
$\geq \gk$.  Thus, we there is a successful coiteration of $(W , \cR )$,
call it $(\cIbar , \cI)$.

It follows from the fact that 
$\cR$ is $\gk$-sound and $(\gk^+)^\cR = (\gk^+)^W$,
that $\cIbar$ is a trivial iteration.
Thus, the coiteration has length OR,
and $W$ the lined up model of $\cI$,
by which we mean that 
$$\Edot^W 
= \bigcup \left\lbrace 
\Edot^{\cM^\cS_\eta} \res \lh (E^\cS_\eta ) \mid \ga \in \OR \right\rbrace$$
Because $W$ is a lower-part $x$-weasel,
it follows that $\cI$ is thorough.
In particular,
$\cI$ has the simple form described in Lemma~2.10.

Let $\ga$ be an ordinal such that 
the drops of $\cI$ are bounded in $\ga$.
Let $\cI'$ be the tail of $\cI$ beginning from stage $\ga$,
and let $\cR'$ be the starting model of $\cI'$.
Thus, $\cI'$ is a thorough iteration of 
$$\cR' = \cM^{\cI'}_0 = \cM^\cI_\ga$$
and $\cI$ has no drops.
Let $$C = \{ \crit (E^{\cI'}_\gb ) \mid \gb \in \OR \}.$$
Then $C$ is a club class of indiscernibles for $W$.
All sufficiently large cardinals are limit points of $C$.
Let $\gm \in \lim (C)$ be a regular cardinal
and
$$U = \left\lbrace A \in W \mid A \subseteq \gm \ \land \ 
\sup ( A \cap C ) = \gm \right\rbrace .$$
Then $\gm > \gk > \gd$ and $U$ is a $\gm$-complete ultrafilter over $W$.
Note that $U$ is just the measure whose trivial completion
is the extender $E^{\cI'}_\gm$,
so that $U$ is normal and weakly amenable to $W$.
\end{proof}

We now briefly discuss how Theorem~2.4 comes up in 
the core model induction of Woodin \cite{W}.
For each rank $V_\ga$ there is a certain lower-part $x_\ga$-weasel $W_\ga$
such that $V_\ga \subset W_\ga$.
Without giving the definition of $W_\ga$, let us say that $W_\ga$ is
the closure of $V_\ga$ under a certain inner model theory operator
corresponding to the least
level of $L ( \mathbb R )$ for which determinacy fails.
As a result, $W_\ga \cap HC$ is $L(\mathbb R )$-definable and $W_\ga$ is 
forcing absolute.
\footnote{Some examples of such weasel 
are given after the statement of Theorem~5.3.}
If there is an $\ga$ such that
$W_\ga$ computes successors of countably closed,
singular cardinals correctly,
then $W_\ga$ can reasonably be construed as the ``true'' core model.
Suppose, on the other hand, that there is no such $\ga$.
Then, Theorem~2.4 implies that for every $\ga$,
there exists $\gm_\ga > \ga$ and a countably complete ultrafilter 
$U_\ga$ on $\cP (\gm_\ga ) \cap W_\ga$ such that $U_\ga$ is weakly amenable
to $W_\ga$.
These properties of $U_\ga$ are
sufficient to develop the Steel core model theory \cite{cmip}
in $W_\ga$ up to $\gm_\ga$.
That is, we have a reasonable interpretation of 
$(K \parallel \gm_\ga)^{W_\ga}$.
The failure of determinacy and the closure of $W_\ga$ under lower parts
implies that $(K \parallel \gm_\ga)^{W_\ga}$ does not reach a Woodin
cardinal.  In this case, it is
$$\bigcup_{\ga \in OR} (K \parallel \ga )^{W_\ga}$$
that is taken as the ``true'' core model.
Because $V_\ga \subset W_\ga$ for all $\ga$,
we can use \cite{MiSchSt} to see that this ``true'' 
core model computes successors of
countably closed, singular cardinals correctly.

Woodin \cite{W} draws on \cite{scales},
\cite{MaSt}, \cite{KW},
\cite{cmip}, \cite{MiSchSt},
and Woodin's theorems on the close relationship between inner 
models with Woodin cardinals and the determinacy of games in $L(\mathbb R)$.

It is possible to use Theorem~3.2 (below) and \cite{MiSch}
instead of Theorem~2.4 and \cite{MiSchSt}, in the core model
induction discussed above.
This gives a core model that computes all successors of
singular cardinals correctly,
based on a failure of $L(\mathbb R )$-determinacy.
\section{Successor cardinals without countable closure}

In this section, we eliminate the hypothesis that $\card ( \gk )$ is
countably closed from Theorem~2.4, and draw a somewhat weaker conclusion.
We shall need a slightly stronger form of iterability.

\begin{definition}\label{3.1}
An $x$-premouse $\cM$ is {\bf $\star$-iterable} 
iff $\cM$ satisfies the following two conditions.
\begin{list}{}{}
\item[(a)]{For every ordinal $\ga$,
$\cM$ is $(\go , \ga )$-iterable in the sense of 
\cite[Definition~2.8]{cmip}.
That is, for every $\ga$,
player II has a winning strategy in the almost normal iteration
game $\cG^* ( \cM , (\go , \ga ))$.}
\item[(b)]{There is no infinite composition 
of thorough iteration trees on $\cM$ 
such that each iteration tree involves a drop.}
\end{list}
\end{definition}

The game $\cG^* ( \cM , (\go , \ga ))$ mentioned in 3.1(a)
is the variation of the usual
iteration game on $\cM$ which allows finitely many violations of the normality
condition on iteration trees.
Another way to express 3.1(b) is to say that there is no infinite sequence 
$\cS_0 {\ }^\frown \cS_1 {\ }^\frown \cdots$
of thorough iteration trees such that 
$\cM^{\cS_0}_0 = \cM,$
and $\cS_n$ has successor length and 
$$\cM^{\cS_{n+1}}_0 \pis \cM_\infty ( \cS_n )$$
for all $n < \go$. 
We remind the reader that thorough iteration trees are
linear, and more, as was established in Lemma~2.1.

\begin{theorem}\label{3.2}
Suppose that $W = L[\Evec , x]$ is an $\star$-iterable,
lower-part $x$-weasel.
Let $\gk$ be a cardinal of $W$ such that 
$$(\cup x)^+ \cup \ha_2 \leq \card (\gk )$$
Put $\gl = ( \gk^+ )^W$ and suppose that $\cf (\gl ) < \card (\gk )$.
Then, given $$\gd < \card (\gk ),$$ there exists $\gm \geq \gd$ and
a normal ultrafilter $U$ on $\cP (\gm ) \cap W$ such
that $U$ is weakly amenable to $W$ and
$\ult (W , U )$ is wellfounded.
\footnote{From the proof of Theorem~3.2, one sees that $U$ can be chosen
so that either $U$ is countably complete,
or $U$ is the measure derived from an inverse transitive collapse
$\gp : \Wbar \lra W$.
In the later case,
$\ult (\ult ( \Wbar , U ) , U )$ is wellfounded, since it embeds into 
$\ult (W , U)$. This is useful in developing the Steel core model theory
\cite{cmip} within $\Wbar$.}
\end{theorem}

Before giving the proof of Theorem~3.2, we isolate a useful
property from the proof of Theorem~2.4.

\begin{definition}\label{3.3}
$\Phi (\cM, \cN)$ holds iff $\cM$ is a lower-part
$x$-premouse, $\cN$ is an $x$-premouse, and there is a successful
coiteration $(\cIbar , \cI)$ of $(\cM , \cN)$ such that $\cIbar$ is trivial
and $\cI$ is thorough.
\end{definition}

Notice that $\Phi (\cM , \cN )$ holds iff $\cM$ is a lower-part $x$-premouse, 
$\cN$ is an $x$-premouse, and for every coiteration
$(\cIbar , \cI)$ of $(\cM , \cN)$,
if $\cIbar$ is trivial and $\cI$ is thorough,
then
$$\lim \left\langle \cM^\cI_\gb \mid \gb < \lh (\cI ) \right\rangle$$ 
is illfounded.
Thus $\Phi (\cM , \cN )$ is highly absolute.

\begin{proof}[Proof of Theorem~3.2]
To do without countable closure, we shall use the same method as
\cite{MiSch}.

Let $\gO > \gl$ be a regular cardinal and
$$\ge = \left( (\cup x) \cup \cf (\gl ) \cup \gd \right)^+ .$$
Then $\ha_2 \leq \ge \leq \card (\gk )$.
Let $\langle X_i \mid i < \ge \rangle$ be a
continuous chain of elementary substructures of $V_{\gO + 1}$ 
such that for all $j < \ge$,
$$(\card (X_j ))^+ = \ge ,$$ 
$$\langle X_i \mid i \leq j \rangle \in X_{j+1},$$ 
and
$$X_j \cap \ge \in \ge .$$
We also have or assume that
$$(\cup x) \cup \cf (\gl ) \cup \gd \cup \{ \gl , \cJ^W_\gO \} \subset X_0 .$$
For $i < \ge$, let $\ge_i = X_i \cap \ge_i$.
Note that $\langle \ge_i \mid i < \ge \rangle$ is a normal sequence
converging to $\ge$.  For $i < \ge$, let $\gp_i : N_i \lra V_{\gO +1}$
be the uncollapse of $X_i$.  So $\crit ( \gp_i ) = \ge_i$.
Say
$\gp_i ( \gk_i ) = \gk$, $\gp ( \gl_i ) = \gl$, and 
$\gp_i ( W_i ) = \cJ^W_\gO$ whenever $i < \ge$.

\bigskip

\noindent
{\bf Case~1.}
{\it There is a stationary set
$$S_0  \subseteq
\left\{ i < \ge \mid \cf ( i ) > \go  \ \land \ i = \ge_i \right\}$$
such that
$i \in S_0$ implies $(\ge_i^+)^{W_i} < (\ge_i^+)^W$.}

\bigskip

In this case, we shall find a arbitrarily large ordinals $\gm$ and 
corresponding ultrafilters $U$ as in the conclusion of the
theorem. Moreover, the ultrafilters will be countably complete.

For $i \in S_0$,
let $(\cTbar_i , \cT_i )$ be the successful coiteration of 
$(W_i , \cJ^W_\gO)$,
and say that this coiteration has length $\theta_i$.
As in the proof of Theorem~2.4,
$\cTbar_i$ is trivial, while $\cT_i$ is dropping and thorough,
whenever $i\in S_0$.
In particular, $\Phi ( W_i , \cJ^W_\gO )$ holds for all $i \in S_0$.

For $i \in S_0$, let $\eta_i$ be the least ordinal $\eta$ such that
$\gn^{\cT_i}_\eta > \gk_i$, if such an ordinal exists, and let $\eta_i$ be
$\theta_i$ otherwise.  And, also, let $\cP_i$ be the longest initial segment of
$\cM^{\cT_i}_{\eta_i}$ (the $\eta_i$'th model of $\cT_i$)
with exactly those subsets of $\gk_i$ that are in
$W_i$.  Then, $\cP_i$ is a $\gk_i$-sound $x$-premouse whenever 
$i \in S_0$.

For $i \in S_0$, let $\cR_i = \ult ( \cP_i , \cE ( \gp_i ) \res \gk )$.
At the corresponding point in the proof of Theorem~2.4, we used countable
completeness to see that the corresponding ultrapower was wellfounded and
iterable.  But here, countable completeness is not available.
Consider an arbitrary $i\in S_0$.
If $\cR_i$ is wellfounded, then it is a $\gk$-sound $x$-premouse that agrees
with $W$ below $\gl$.  If, in addition, $\cR_i$ is iterable above $\gk$,
then $\Phi ( W , \cR_i )$ holds, giving us the indiscernibles for $W$,
and we are done as in Theorem~2.4.
Thus, we may assume that $\Phi (W , \cR_i )$ fails for all $i\in S_0$.

Let $S$ be a stationary subset of 
$S_0$ on which the choice functions
$i \mapsto \eta_i$
and 
$i \mapsto n(\cP_i , \gk_i )$
are constant.  Such a set $S$ exists by Fodor's lemma \cite[Lemma~1.1]{MiSch}.
Say $\eta_i = \eta$ and $ n(\cP_i , \gk_i ) = n$ whenever $i \in S$.

So if $i \in S$ and $\cR_i$ is wellfounded, then $n (\cR_i , \gk )=n$,
$\cR_i$ agrees with $W$ below $\gl$, but $\cR_i$ is not iterable above
$\gk$.  Even if
$\cR_i$ is illfounded, $\gl$ is in the wellfounded part of $\cR_i$ and
$\cR_i$ agrees with $W$ below $\gl$. (We remark that, without loss of
generality, either $\cR_i$ is not wellfounded for all $i \in S$, or $\cR_i$
is wellfounded but not iterable for all $i \in S$.)

For the rest of the proof of the theorem under Case~1, 
fix 
$$j \in \lim (S ) \cap S.$$
Let $\psi : M \lra V_{\gO + 2}$ be elementary with
$M$ countable and transitive, and everything relevant in the range of $\psi$.
Say $\psi (x' ) = x$, $\psi ( \gO' ) = \gO$, $\psi (W' ) = \cJ^W_\gO$, 
$\psi ( \cP' ) = \cP_j$,
$\psi ( \cR' ) = \cR_j$,
and $\psi (\gk' ) =\gk$.
Then $\Phi (W' , \cR')$ fails in $M$,
hence in $V$, so $\cR'$ is not iterable above $\gk'$
(possibly because $\cR'$ is not wellfounded).

\bigskip

\noindent
{\bf Claim~3.2.1.}
{\it $\ran (\psi ) \cap \gl$ is unbounded in $\gl$.}

\bigskip

\begin{proof}[Sketch of Claim~3.2.1]
Suppose to the contrary that 
$\ran ( \psi ) \cap \gl$ is bounded in $\gl$.
Let $\gptil : \cP_j \lra \cR_j$ be the ultrapower map.
Then $\gptil'' \gl_j$ is unbounded in $\gl$ and $\gptil \in \ran (\psi )$.
Moreover, the following diagram commutes.\\
$$\begin{diagram}
\node{\cP_j}
\arrow[3]{e,t}{\gptil}
\node{}
\node{}
\node{\cR_j}\\
\node{\cP'}
\arrow{n,l}{\psi}
\arrow[3]{e,b}{\psi^{-1}(\gptil )}
\node{}
\node{}
\node{\cR'}
\arrow{n,r}{\psi}\\
\end{diagram}$$
Therefore $\ran (\psi ) \cap \gl_j$ is bounded in $\gl_j$.

We shall use \cite[\S2]{square} to find a proper initial segment
$\cL$ of $\cP_j$ and an almost $\gS_{n+1}$-embedding $\varphi : \cL \lra
\cP_j$ such that 
$$\ran (\psi ) \cap | \cP_j | \subset \ran (\varphi ),$$
$\varphi \in | \cP_j |$, and 
$$\ran (\psi ) \cap |\cR_j | \subset \ran
(\gptil (\varphi )).$$  
Moreover, we shall have $\gr_{n +1} (\cL ) = \gk_j$
so that $\cL \pis \cP_j \parallel \gl_j$ and 
$$\gptil ( \cL ) \pis \cR_j \parallel \gl.$$
Then
$\gp_j ( \cL ) = \gptil ( \cL )$
by the agreement of
$\gp_j$ and $\gptil$.
So $\gp_j (\cL ) \pis W$
by the agreement of $\cR_j$ and $W$, 
Therefore
$\gp_j ( \cL )$ is iterable.

We shall have the following diagram.\\
$$\begin{diagram}
\node{\cP_j}
\arrow[3]{e,t}{\gptil}
\node{}
\node{}
\node{\cR_j}\\
\node{}
\node{\cL}
\arrow{nw,t}{\varphi}
\arrow{e,t}{\gp_j}
\node{\gp_j ( \cL )}
\arrow{ne,t}{\gptil (\varphi)}
\node{}\\
\node{\cP'}
\arrow[2]{n,l}{\psi}
\arrow{ne}{}
\arrow[3]{e,b}{\psi^{-1}(\gptil )}
\node{}
\node{}
\node{\cR'}
\arrow{nw,b}{\psitil}
\arrow[2]{n,r}{\psi}\\
\end{diagram}$$

We have labeled 
$(\gptil (\varphi ))^{-1} \circ \psi $ as $\psitil $ in the
diagram above.  The map $\psitil$ is a weak $n$-embedding which is
$\gS_{n+1}$-elementary from $\cR'$ into $\gp_j (\cL )$.  
Because $\gp_j (\cL )$ is iterable,
the existence of a map such as $\psi$ implies that $\cR'$ is iterable. 
But $\cR'$ is not iterable above $\gk'$, 
while $\gp_j (\cL )$ is iterable.
This is the desired contradiction.

It remains to construct $\varphi: \cL \lra \cP_j $ as described above.  
Let $\gr = \gr_n ( \cP_j )$.
For notational simplicity alone,
let us restrict our attention to the case $n=0$ and $\cP_j$ is not type III,
so that $\gr = \OR^{\cP_j}$
(the general case can be handled similarly
by working over the $\gS_n$ coding structure of $\cP_j$, as is usual in Jensen
fine structure).
In the proof of \cite[Theorem~2.9]{square}, 
it is shown how to obtain a parameter $r$ in $\cP_j$ such that 
$p_1 ( \cP_j )$ is included in $r$, with the property that
for all $\ga < \gk_j$, if we set
$$\grbar ( \ga ) = \sup \left( 
H^{\cP_j}_1 \left( \ga \cup r \right) \cap \gr \right)$$
and
$$\glbar (\ga ) = 
\sup \left( H^{\cP_j}_1 \left( \ga \cup r \right) \cap \gl_j \right) ,$$
then
$$\grbar ( \ga ) < \gr \ \Llra \ \glbar ( \ga ) < \gl_j \ .$$
Without repeating the details, let us remark that
the main facts used in defining $r$ are the uncountable cofinality of 
$\gl_j$ and the $\gk_j$-soundness of $\cP_j$.

For any $\ga < \gk_j$, let
$$\cL_\ga = \cH^{\cP_j \res \grbar ( \ga )}_1 ( \glbar (\ga ) \cup r )$$
and
$$\varphi_\ga : \cL_\ga \lra \cP_j$$ 
be the inverse transitive collapse.

In what follows, we consider mainly $\ga < \gk_j$ for which 
$\cL_\ga \not= \cP_j$.  Equivalently, $\grbar ( \ga ) \not= \gr$.
Equivalently, $\glbar ( \ga ) \not= \gl_j$. 
Call such $\ga$ {\bf reasonable}.
If $\ga$ is reasonable, then $\varphi_\ga \in |\cP_j |$.

We claim that for all sufficiently large reasonable $\ga < \gk_j$,
setting $\varphi = \varphi_\ga$ and $\cL = \cL_\ga$ works.

We first argue that for all $\ga < \gk_j$,
$\cL_\ga$ is an $x$-premouse.
In the notation of \cite{square},
$\varphi_\ga$ is an almost $\gS_1$-embedding.
General facts about such embeddings reduce the problem of seeing that
$\cL_\ga$ is an $x$-premouse to seeing that,
if $\cP_j$ is active, then
$\Fdot^{\cL_\ga}$ measures all sets in $\cL_\ga$.
In the terminology of \cite{MiSchSt}, we must see that $\cL_\ga$ is
an $x$-premouse, and not merely an $x$-protomouse.
Note that this is a real concern,
given that if $\cP_j$ is active and $\gs < \OR^{\cP_j}$,
then $\cP_j \res \gs$ is certainly not an $x$-premouse.

Suppose that $\cP_j$ is active.  In order to see that $\Fdot^{\cL_\ga}$
measures all sets in $\cL_\ga$, it is enough to see that 
$\gmdot^{\cP_j} \geq \gk_j$.
Suppose to the contrary that $\gmdot^{\cP_j} < \gk_j$.
Since $\gl_j$ is a cardinal in $\cP_j$,
by the initial segment condition on $\cP_j$,
there is a proper initial segment $G$ of $\Fdot^{\cP_j}$ which is on the
$\cP_j$-sequence below $\gl_j$.
Recall that $\cP_j$ agrees with $W_j$ below $\gl_j$
and that $\gl_j$ is the successor cardinal of $\gk_j$ in both 
$\cP_j$ and $W_j$.  Therefore, $G$ is an extender on the $W_j$-sequence
which measures all sets in $W_j$.  This is a contradiction because $W_j$
is a lower-part $x$-weasel.

We next claim that $\cL_\ga \pis \cP_j$
for sufficiently large reasonable $\ga < \gk_j$.
\footnote{The analogous fact was overlooked in \cite{square};
a slightly simpler proof of \cite[Theorem~2.9]{square} is possible.}
First, 
since $r$ contains $p_1 ( \cP_j )$ and $\cP_j$ is $\gk_j$-sound,
for all sufficiently large reasonable $\ga < \gk$,
$$\crit ( \varphi_\ga ) = \glbar (\ga ),$$
$$\cL_\ga \mbox{ is $1$-sound,}$$
$$\gr_1 ( \cL_\ga ) = \gk_j ,$$
$$\varphi_\ga ( p_1 ( \cL_\ga )) = p_1 ( \cP_j ),$$
and
$$p_1 ( \cL_\ga ) = (\varphi_\ga )^{-1} \left( p_1 ( \cP_j )\right) .$$
It enough to take $\ga$ large enough so that $r$
and the solidity witnesses
for $p_1 (\cP_j )$ are all elements of 
$$H^{\cP_j \res \grbar ( \ga )}_1 ( \gk_j \cup p_1 ( \cP_j ) ).$$

If $\cP_j$ were $1$-small, then we could apply the condensation result
\cite[Theorem~2.6]{square} to see that 
$\cL_\ga \pis \cP_j$
for all sufficiently large reasonable $\ga < \gk_j$.
But there
is no smallness condition on $W$, hence not on $\cP_j$, so a small argument is
required.

The fact that $W_j$ is a lower-part $x$-weasel that agrees with $\cP_j$
below $\gl_j$ implies that in any coiteration of $(\cL_\ga , \cP_j )$,
no extender with critical point $< \gk_j$ is used.  Since $\cP_j$ is
iterable above $\gk_j$, and the embedding $\varphi_\ga : \cL_\ga \lra \cP_j$
is sufficiently elementary, $\cL_\ga$ is iterable above $\gk_j$.
So there is a successful coiteration $(\cSbar , \cS )$ of $(\cL_\ga , \cP_j )$.
But now the $\gk_j$-soundness of both $\cP_j$ and $\cL_\ga$ can be used to
see that $\cL_\ga \pis \cP_j$.

Now fix $\cL = \cL_\ga$ and $\varphi = \varphi_\ga$ for some 
reasonable $\ga < \gk_j$ which is large enough so that
$\cL_\ga \pis \cP_j$ and 
$$\sup \left( \ran (\psi ) \cap \gr \right) < \gpbar ( \ga ).$$
This last condition guarantees that 
$\ran (\psi ) \cap |\cP_j | \subset \ran (\varphi )$
and 
$$\ran (\psi ) \cap |\cR_j | \subset \ran (\gptil ( \varphi  ) ).$$
All the requirements of $\varphi$ and $\cL$ set out at the beginning of the
proof of Claim~3.2.1 are now plain, and so we are done with that claim.
\end{proof}

We continue the proof of Theorem~3.2.  Since $\psi$ is countable,
there is an $i \in j \cap S$ such that $\ran (\psi ) \cap X_j \subset X_i$.
Here we are using that $j \in \lim (S)$ and $\cf ( j ) > \go$.
Fix such an $i$ for the rest of the proof of the theorem under Case~1.
Let $\gp_{i,j} : N_i \lra N_j$ be the uncollapse of $(\gp_j )^{-1}(X_i )$.
Since $\ran ( \psi ) \cap \gl$ 
is unbounded in $\gl$ (by Claim~3.2.1),
we have that
$$\sup \left( (\gp_{i,j} ) '' \gl_i \right) = \gp_{i,j}
(\gl_i ) = \gl_j \ . $$
Arguing as in \cite[Lemma~2.1]{MiSch}, we find an 
$x$-premouse $\cP_j^*$ such that $\cP_j^*$ agrees with $W_i$ below $\gl_i$
and $$\cP_j = \ult ( \cP_j^* , \cE (\gp_{i,j}) \res \gk_j )$$
Moreover,
$((\gk_i)^+ )^{\cP_j^*} = \gl_i$
and $\cP_j^*$ is $\gk_i$-sound with $n ( \cP_j^* , \gk_i ) = n$.
In \cite{MiSch},
$\cP^*_j$ is what is called the ``pull-back'' of $\cP_j$ by $\gp_{i,j}$.

We claim that $\cP_i = \cP_j^*$.  Note that any coiteration of $\cP_i$ and
$\cP_j^*$ involves only extenders with critical points at least $\gk_i$.
This follows from the initial segment condition,
using the fact that both $\cP_i$ and $\cP_j^*$ are premice which
agree with $W_i$ up to $\gl_i$,
$W_i$ is a lower-part $x$-premouse, and $\gl_i$ is the cardinal successor of
$\gk_i$ in all three premice.
Since both $\cP_i$ and $\cP_j^*$ are iterable above $\gk_i$,
there is a successful coiteration $(\cSbar , \cS)$ of $(\cP_i , \cP_j^* )$.
But now the standard argument using the $\gk_i$-soundness of $\cP_i$ and
$\cP_j^*$ shows that the coiteration is trivial, and hence $\cP_i =
\cP_j^*$.

Recall that $\cR_i = \ult ( W_i , \cE (\gp_j) \res \gk)$ and that
$\Phi ( W , \cR_i )$ fails.
In fact, $\Phi ( \cJ^W_\gO , \cR_i )$ fails.
So 
$$N_j \models 
\neg \Phi \left( \left(\gp_j \right)^{-1} \left( \cJ^W_\gO \right) , 
\left(\gp_j \right)^{-1} \left( \ult \left( W_i , \cE (\gp_i) \res \gk \right)
\right) \right).$$
But $\gp_j ( W_j ) = \cJ^W_\gO$,
$$\gp_j \left( \gp_{i,j} \res \gk_j \right) = \gp_i \res \gk ,$$
and, as we have seen,
$$\cP_j = \ult ( \cP_i , \cE (\gp_{i,j}) \res \gk_j ).$$
Therefore
$$N_j \models \neg\Phi ( W_j , \cP_j ).$$
By upward absoluteness, $\Phi ( W_j , \cP_j )$ fails.

But notice that $\Phi ( W_j , \cP_j )$ holds trivially,
since $P_j$ is an initial segment of a model on $\cT_j$ and
$(\cTbar_j , \cT_j)$ is a coiteration of $(W_j , W)$.
Thus, we have a contradiction, and the proof of Theorem~3.2 under Case~1
is complete.

\bigskip

\noindent
{\bf Case~2.}
{\it There is a stationary set
$$S \subseteq 
\left\{ i < \ge \mid \cf ( i ) > \go  \ \land \ i = \ge_i \right\}$$
such that
$i \in S$ implies $(\ge_i^+)^{W_i} = (\ge_i^+)^W$.}

\bigskip

In this case, we shall find an ordinal $\gm$ with $\gd < \gm < \card ( \gk )$,
and an ultrafilter $U$ on $\cP ( \gm ) \cap W$ such that $U$ is weakly
amenable to $W$, and $\gO$ contained in the wellfounded part of 
$\ult ( W,U)$.  Since $\gO$ can be taken arbitrarily large, and there are
only set many possible ultrafilters $U$, there will be at least one $U$ that 
satisfies the conclusion of the theorem.

As in the proof of Theorem~2.4, $W$ and $W_i$ have the same subsets of
$\ge_i$ whenever $i \in S$.
Let $\cP = \cJ^W_\gO$.  Note that $\cP$ is a lower-part $x$-premouse,
since $\gO$ is a cardinal of $W$.
Let $\cR_i = \ult (W_i , \cE (\gp_i) \res \ge )$ whenever $i \in S$.
In other words, $\cR_i$ is the ultrapower of $\cP$ by the superstrong
extender derived from $\gp_i$.

It is enough to see that there exists an $i \in S$ such that $\cR_i$ is
wellfounded, for the following reasons.
Let $U_i$ be the measure derived from $\gp_i$, that is,
$$U_i = \left\lbrace x \in W_i \mid x \subseteq \ge_i \ \land \ 
\ge_i \in \gp_i ( x ) \right\rbrace.$$
Then $U_i$ is an ultrafilter over $\cP ( \ge_i ) \cap W$
and $U_i$ is weakly amenable to $W$.
Also, there is an embedding of $\ult ( \cP , U_i )$ into $\cR$.
Therefore, if $\cR_i$ is wellfounded, then so is 
$\ult (\cP ,U_i )$.

Assume towards a contradiction that $\cR_i$ is illfounded whenever 
$i \in S$.

Fix $j \in S \cap \lim ( S)$.  We shall show that $\cR_j$ is wellfounded by
an argument similar to that used in Case~1, letting $\cP$ play the role of
$\cP_i$ for all $i \in S$.

Let $\psi : M \lra V_{\gO + 2}$ be elementary with $M$ countable and
transitive, and everything relevant in the range of $\psi$.
Say $\psi (x' ) = x$, $\psi ( \gO' ) = \gO$, $\psi (W' ) = \cJ^W_\gO$,
$\psi ( \cP' ) = \cP$,
$\psi ( \cR' ) = \cR_j$,
and
$\psi (\gk' ) =\gk$.
Then $\cR'$ is illfounded.

For the rest of the proof of the theorem,
fix $i \in j \cap S$ such that $\ran (\psi ) \cap X_j \subset X_i$.
Let $\gp_{i,j} : N_i \lra N_j$ be the uncollapse of $(\gp_j )^{-1}(X_i )$.
Arguing as in \cite[Lemma~2.2]{MiSch}, we find a ``pull-back'' $\cP^*$ 
of $\cP$ by $\gp_{i,j}$.  That is, we find an
$x$-premouse $\cP^*$ such that $\cP^*$ agrees with $W$ below 
$$\ge^* = ((\ge_i )^+)^{\cP^*},$$
and the following commutative diagram.\\
$$\begin{diagram}
\node{\cP^*}
\arrow{e,t}{}
\arrow{se,b}{}
\node{\cP}
\arrow{e,t}{}
\node{\cR_j}\\
\node{}
\node{\ult ( \cP^* , \cE (\gp_{i,j}) \res \ge_j )}
\arrow{n,l}{}
\arrow{e,b}{}
\node{\ult ( \cP^* , \cE (\gp_i) \res \ge )}
\arrow{n,l}{}
\node{\cR'}
\arrow{nw,t}{}
\arrow{w,b}{}\\
\end{diagram}$$
Thus $\ult ( \cP^* , \cE (\gp_{i,j}) \res \ge_j )$ is wellfounded and
iterable,
while $\ult ( \cP^* , \cE (\gp_i) \res \ge )$ is illfounded.

Recall that $\cR_i$ is illfounded.
Therefore
$$N_j \models \mbox{``$\gp_j^{-1} ( \cR_i )$ is illfounded''}.$$
Hence
$$\gp_j^{-1} ( \cR_i ) 
= \ult (\gp_j^{-1} ( \cP ) , \cE (\gp_{i,j}) \res \ge_j ) $$
is illfounded.
The shift lemma gives us the following diagram.\\
$$\begin{diagram}
\node{\ult (\gp_j^{-1} ( \cP ) , \cE (\gp_{i,j}) \res \ge_j )}
\arrow{e,b}{}
\node{\ult (\cP , \cE (\gp_{i,j}) \res \ge_j )}\\
\node{\gp_j^{-1} ( \cP )}
\arrow{e,b}{\gp_j}
\arrow{n,l}{}
\node{\cP}
\arrow{n,r}{}\\
\end{diagram}$$
Thus both $\cR_i = \ult ( \cP , \cE (\gp_i) \res \ge )$ and
$\ult (\cP , \cE (\gp_{i,j}) \res \ge_j )$ are illfounded.

We define a ``minimal'' $x$-premouse $\cQ$ as follows.
This is where $\star$-iterability is used rather than just iterability.
Let $\cQ_0$ be the least initial segment $\cQ'$ of $\cP$ above $\ge_i$
such that $\cQ'$ is a lower-part $x$-premouse and
$\ult ( \cQ' , \cE (\gp_i) \res \ge )$ is illfounded.
Suppose that $0 < \ell < \go$ and for all $k < \ell$,
we have defined $\cQ_{k+1}$ and a non-trivial, 
thorough iteration tree $\cS_{k+1}$
on $\cQ_k$ such that $\cQ_{k+1} \is \cM_\infty (\cS_{k+1} )$.
If there is a non-trivial thorough iteration tree $\cS$ of successor length on
$\cQ_\ell$ such that all critical points of extenders used on $\cS$ are at
least $\ge_i$,
and there is a proper 
initial segment $\cQ'$ of $\cM_\infty ( \cS )$ such that
$\cQ'$ is a lower-part $x$-premouse and 
$\ult ( \cQ' , \cE (\gp_i) \res \ge )$ is illfounded,
then let $\cS_{\ell + 1}$ and $\cQ_{\ell + 1}$ be the least
such $\cS$ and $\cQ'$.
\footnote{Note that thorough iteration trees on a fixed $x$-premouse 
are indexed by finite sets of ordinals.}
Otherwise, leave $\cS_{\ell + 1}$ and $\cQ_{\ell + 1}$ undefined.
By the $\star$-iterability of $W$ (because there is no infinite composition of
thorough iteration trees on $W$ such that each iteration tree involves a
drop), there is a $k < \go$ such that $\cQ_{k+1}$ is not defined.
Let $\cQ = \cQ_k$.  

We have that $\cQ$ is a lower-part $x$-premouse, and 
$\cQ$ is iterable since $W$ is $\star$-iterable.
Note that $\ult ( \cQ , \cE (\gp_i) \res \ge )$ is illfounded,
while if $\cS$ is any thorough iteration tree on $\cQ$ with critical points
at least $\ge_i$ and $\cQ'$ is any proper
initial segment of $\cM_\infty ( \cS )$,
then $\ult ( \cQ' , \cE (\gp_i) \res \ge )$ is wellfounded.
Also note that $\cQ \in X_j$.
Essentially the same argument that showed
that $$\ult ( \cP , \cE (\gp_{i,j}) \res \ge_j )$$ is illfounded shows that
$\ult ( \cQ , \cE (\gp_{i,j}) \res \ge_j )$ is illfounded.
\footnote{Just replace $\cP$ by $\cQ$ and $\cR_i$ by 
$\ult ( \cQ , \cE (\gp_i) \res \ge)$ two paragraphs above.}

Now let $(\cSbar , \cS)$ be a successful coiteration of $(\cQ , \cP^* )$.
Both $\cQ$ and $\cP^*$ are lower-part $x$-premice, so one of $\cSbar$ and
$\cS$ is trivial, and the other is thorough.
The minimality of $\cQ$ described above, and the fact that 
$\ult (\cP^* , \cE (\gp_i) \res \ge )$ is illfounded guarantees that $\cSbar$
is trivial.
Then $\cQ \is \cM_\infty ( \cS )$.
Consider the following figure.\\
$$\begin{diagram}
\node{}
\node{}
\node{\ult(\cM_\infty ( \cS ), \cE (\gp_{i,j}) \res \ge_j)}
\arrow{se,t}{}
\node{}
\node{}\\
\node{\cQ}
\arrow{e,t}{\is}
\node{\cM_\infty ( \cS )}
\arrow{ne,t}{}
\arrow[2]{e,b}{\gs}
\node{}
\node{\cM_\infty ((\gp^*)''\cS )}\\
\node{}\\
\node{}
\node{\cP^*}
\arrow[2]{n,l,..}{\cS}
\arrow[2]{e,b}{\gp^*}
\node{}
\node{\ult (\cP^* , \cE (\gp_{i,j}) \res \ge_j)}
\arrow[2]{n,r,..}{(\gp^*)''\cS}\\
\end{diagram}$$
By $\gp^*$ we mean the ultrapower map.  The vertical dashed lines
represent iteration trees, but not necessarily embeddings.
Recall that $\ult (\cP^* , \cE (\gp_{i,j}) \res \ge_j)$ is iterable.
We copy $\cS$ using $\gp^*$, and let $\gs$ be the final copying map.
Then $\gs \res \ge_i = \gp \res \ge_i$ since all critical points of extenders
used on $\cS$ are at least $\ge_i$. 
So 
$$\ult(\cM_\infty ( \cS ), \cE (\gp_{i,j}) \res \ge_j) = 
\ult(\cM_\infty ( \cS ), \cE (\gs)  \res \ge_j)$$
embeds into $\cM_\infty ((\gp^*)''\cS )$ as indicated.
In particular,
$$\ult(\cM_\infty ( \cS ), \cE (\gp_{i,j}) \res \ge_j)$$
is wellfounded.

But $\cQ \is \cM_\infty ( \cS )$
and $\ult(\cQ , \cE (\gp_{i,j}) \res \ge_j)$ is illfounded,
so $$\ult(\cM_\infty ( \cS ), \cE (\gp_{i,j}) \res \ge_j)$$ is illfounded.
This contradiction completes our consideration of Case~2,
and the proof of Theorem~3.2.
\end{proof}
\section{A weak iterability condition}

Notice that in \S\S2-3, the iterability of $W$ was used to show that all
the relevant iteration trees are linear, even thorough.
But, the linear iterability of $W$ does not suffice for those arguments.

Nevertheless, the assumption that $W$ 
is $\star$-iterable is more than is required for the proofs
of Theorems~2.4 and 3.2.  In this section, we sketch why internal $\star$-iterability
for countable premice elementarily realizable in a level of $W$ suffices for
the arguments in the preceding sections.

\begin{definition}\label{4.1}
Let $W$ be an $x$-weasel and $\ga$ an ordinal.
We say that a structure $\cP$ is {\bf realizable in $W$} iff
$\cP$ is an $\xdot^\cP$-premouse and there exists an ordinal $\gO$
and an elementary embedding of $\cP$ into $\cJ^W_\gO$.
We say that $W$ is {\bf $\star \star$-iterable} iff the following
conditions hold for every countable $\cP$ which is realizable in $W$.
\begin{list}{}{}
\item[(a)]{For every ordinal $\ga \leq \go_1$,
$\cP$ is $(\go , \ga)$-iterable in the sense of \cite[Definition 2.8]{cmip}.}
\item[(b)]{There is no infinite composition 
of thorough iteration trees on $\cP$ 
such that each iteration tree involves a drop.}
\end{list}
We say that $W$ is {\bf internally $\star \star$-iterable} iff
$$W \models \mbox{``$L[\Edot , \xdot ]$ is $\star \star$-iterable''}.$$
\end{definition}

\begin{theorem}\label{4.2}
Suppose that $W = L[\Evec , x]$ is an internally
$\star \star$-iterable, lower-part $x$-weasel.
Let $\gk$ be a cardinal of $W$ such that 
$$(\cup x )^+ \cup \ha_2 \leq \card (\gk )$$
Put $\gl = ( \gk^+ )^W$ and suppose that $\cf (\gl ) < \card (\gk )$.
Then, given $$\gd < \card (\gk ),$$ there exists $\gm \geq \gd$ and
a normal ultrafilter $U$ on $\cP (\gm ) \cap W$ such
that $U$ is weakly amenable to $W$ and $\ult ( W , U )$ is wellfounded.
\end{theorem}

Theorem~4.2 is just Theorem~3.2 restated 
with $\star$-iterability weakened to internal $\star \star$-iterablity;
we now sketch how to modify the proof of the former to yield the latter.
Let us consider the first time that iterability is used in the proof
of Theorem~3.2, which is to see that in Case~1 of the proof,
$(\cTbar_i , \cT )$ witnesses that $\Phi ( W_i , \cJ^W_\gO )$
holds whenever $i \in S_0$.  Fix $i \in S_0$ and suppose that
$\Phi ( W_i , \cJ^W_\gO )$ fails.
By the absoluteness of wellfoundedness, the objects listed below exist in $W$,
because the objects exist in $V$ and there is a tree in $W$ that searches
for them.
\begin{itemize}
\item{a countable, transitive model $M$ of a large chunk of ZFC}
\item{a set of ordinals $x^M$ in $M$, a lower-part $x^M$-premouse 
$W^M$ in $M$, and ordinals $\gd^M < \gk^M < \gl^M < \gO^M$ in $M$, 
which are related as in the hypothesis of Theorem~2.4 from the point of view
of $M$}
\item{an internally approachable chain of 
elementary substructures of $(V_{\gO^M + 1} )^M$ of the sort that
arises when we attempt to run the proof of Theorem~2.4 inside $M$}
\item{a map $\gp^M : N^M \lra (V_{\gO^M + 1} )^M$ in $M$ which comes from
collapsing
one of the substructures on the internally approachable chain of $M$ 
to be transitive}
\item{an $x^M$-premouse $\Wbar^M$ such that $\gp^M ( \Wbar^M ) = W^M$
and 
$$M \models \neg \Phi ( \Wbar^M , W^M )$$}
\item{an elementary embedding $\gt : W^M \lra \cJ^W_\gO$}
\end{itemize}

The existence of elementary embeddings
$\gt \circ \gp : \Wbar^M \lra \cJ^W_\gO$ and 
$\gt : W^M \lra \cJ^W_\gO$ 
in $W$
implies 
that $\Wbar^M$ and $W^M$ are both countable, realizable $x^M$-premice in
$W$, and are therefore $\star$-iterable in $W$.
In particular, they are iterable in $W$, so there is a 
successful coiteration $(\cTbar , \cT )$ of $(\Wbar^M , W^M )$ in $W$.
Arguing as in the proof of Theorem~2.4, we see that $\cTbar$ is trivial and
that $\cT$ is thorough (in $W$, hence in $V$).
But then $(\cTbar , \cT)$ witnesses that 
$M \models \Phi ( \Wbar^M , W^M )$, which is a contradiction.

There are several other 
applications of iterability in the proof of Theorem~3.2.
All can be reduced to applications of internal 
$\star \star$-iterability using the
preceding idea.  We leave the details to the reader.
\section{The Jensen covering property}

\begin{definition}\label{5.1}
If $W$ is an $x$-weasel and $\gk$ is an ordinal, then CP$(W, \gk )$ holds 
iff $\cup x < \card ( \gk )$ and 
for every uncountable $A \subseteq \gk$ with
$\cup x \subseteq A$, there exists $B \in W$ such that 
$A \subseteq B$ and $\card (B) = \card (A)$.
\end{definition}

\begin{definition}\label{5.2}
An $x$-weasel $W$ satisfies {\bf strong condensation} 
iff for every $x$-premouse $\Wbar$,
if there is an elementary embedding from $\Wbar$ into a rank
initial segment of $W$, then $\Wbar \pis W$.
\end{definition}

\begin{theorem}\label{5.3}
Let $W = L [ \Evec , x]$ 
be an internally $\star \star $-iterable, lower-part $x$-weasel.
Suppose that $\gk$ is the least ordinal such that 
CP$(W, \gk)$ fails.  Let $A$ be a witness to the failure of CP$(W , \gk )$.
Then one of the following holds.
\begin{list}{}{}
\item[(a)]{There is an ordinal $\gm > \card (A)$ and a
normal ultrafilter $U$ on $\cP ( \gm ) \cap W$ such that $U$ is weakly
amenable to $W$ and $\ult ( W , U )$ is wellfounded.}
\item[(b)]{There exists $\cR \pis W$ with $\gr_\go ( \cR ) = \gk$,
and there exists a set $C \subset \gk$ of indiscernibles for 
$W \parallel \gk$ such that $C$ is unbounded in $\gk$
and the last extender of $\cR$ is a measure generated by $C$. 
In addition, $$A \subset H^\cR_\go ( \gkbar \cup C )$$
for some $\gkbar < \gk$.}
\end{list}
Moreover, if $W$ satisfies strong condensation,
then it is (a) that holds.
\end{theorem}

The simplest weasel that satisfies strong condensation is $L$;
thus Theorem~5.3 extends Jensen's covering lemma for $L$.
A more complicated example of a lower-part 
weasel that satisfies strong condensation
is the minimal model closed under sharps for sets (assuming sharps exist).
Another example is the minimal model closed under the operator 
which sends a set $X$ to $M_1 (X)^\#$, 
the sharp for the iterable one Woodin cardinal model built over $X$
(assuming the operator is total).
These examples of $x$-weasels are typical of those that come up in \cite{W},
in that they are minimal with respect to a certain lower-part inner model
operator, and hence satisfy strong condensation.

We remark that 
there are examples showing the need for clause~(b) in Theorem~5.3.

\begin{proof}[Proof of Theorem 5.3]
Let us assume that $W$ is 
$\star$-iterable; the reduction to internal 
$\star \star $-iterability is as outlined in \S4.

Let $\gk$ be least such that CP$(W,\gk )$ fails,
and let $A \subset \gk$ witness the failure of CP$(W , \gk )$.
Let $\ge = (\card (A))^+$.
Then
$$\cf ( \gk ) \leq \card (A) < \ge \leq \card ( \gk ) \leq \gk = \sup ( A ).$$
Since $\gk$ is least such that CP$(W , \gk )$ fails,
$\gk$ is a cardinal in $W$.

In this proof, $\gk$ will play the role that $\gl$ played in
the proofs of our other theorems.
Note that by Theorem 3.2, we may assume that
$\gk$ is a limit cardinal of $W$.

Let $\gO > \gk$ be a regular cardinal.
Let $\langle X_i \mid i < \ge \rangle$ be a continuous chain of elementary
substructures of $V_{\gO+1}$ such that for all $j < \ge$,
$$\card (X_j ) = \card (A),$$
$$\langle X_i \mid i \leq j \rangle \in X_{j+1} ,$$
and
$$X_j \cap \ge \in \ge .$$
For $i < \ge$, let $\ge_i = X_i \cap \ge_i$ and $\gp_i : N_i \lra V_{\gO+1}$
be the uncollapse of $X_i$, so that $\crit ( \gp_i ) =\ge_i$.
We also have or assume that
$$(\cup x ) \cup \card (A) \cup \{ \gk , \cJ^W_\gO , A \} \in X_0.$$ 
Say $\gp_i ( \gk_i ) = \gk$ and $\gp_i ( W_i ) =
\cJ^W_\gO$ whenever $i < \ge$.

\bigskip

\noindent
{\bf Case~1.}
{\it There is a stationary set
$$S_0 \subseteq 
\left\{ i < \ge \mid \cf (i ) > \go \ \land \ i=\ge_i \right\}$$
such that $i \in S_0 $ implies $(\ge_i^+ )^{W_i} < (\ge_i^+ )^W$.}

\bigskip

For $i \in S_0$, let $(\cTbar_i , \cT_i)$ be a successful coiteration of
$(W_i , \cJ^W_\gO )$, and say this coiteration has length $\theta_i$.  
Then $\cTbar_i$ is trivial and $\cT_i$ is thorough
whenever $i \in S_0$.  In particular, $\Phi (W_i , \cJ^W_\gO )$ holds.

Notice that if $W$ satisfies strong condensation, then $\Wbar \pis W$,
so $\cT_i$ is also trivial.

For $i \in S_0$, let $\eta_i$ be the least ordinal $\eta$ such that $\crit
(E^\cT_\eta ) \geq \gk_i$, if such an ordinal exists, and let $\eta_i$ be
$\theta_i$ otherwise.  For $i \in S_0$, 
let $\cP_i$ be the longest initial segment of
$\cM^{\cT_i}_{\eta_i}$ (the $\eta_i$'th model of $\cT_i$)
with precisely the same bounded subsets of $\gk_i$ as
$W_i$.
Then each $\cP_i$ is a
$\gk_i$-sound $x$-premouse.

For $i \in S_0$, let $\cR_i = \ult (\cP_i , \cE ( \gp_i ) \res \gk )$.
Recall that this ultrapower is defined using coordinates 
$a \in [\gk ]^{< \go}$,
and consists of equivalence classes of certain $\cP_i$-definable functions 
of the form $f : [ \ga ]^{|a|} \lra |\cP_i |$ with $\ga < \gk_i$.
This makes sense because $\cP_i$ and $W_i$ have the same bounded
subsets of $\gk_i$, and because $\gp_i$ is continuous at $\gk_i$.

\bigskip 

\noindent
{\bf Claim 5.3.1.}  
{\it There is a $i \in S_0$ such that
either $\Phi (\cJ^W_\gO , \cR_i)$ holds or $\cR_i \pis W $.}

\bigskip

Before proving Claim~5.3.1, let us see how to use it to finish the proof of
Theorem~5.3 under Case~1.  Fix $i$ as in Claim 5.3.1.
To ease the notation,
let us put $\gp = \gp_i$, $\Wbar = W_i$, $\gkbar = \gk_i$,
$\eta = \eta_i$, $\cT = \cT_i$, $\cP = \cP_i$, and $\cR = \cR_i$
until we start the proof of Claim~5.3.1.

If $\Phi ( W , \cR )$ holds, then (a) holds by an argument as in the proof
of Theorem~2.4.  We can also conclude that (a) holds if
for a proper class of $\gO$, there is a corresponding
$i ( \gO )$ such that $\Phi (\cJ^W_{\gO} , \cR_{i ( \gO )})$ holds.
This is by a counting argument, using that there are only set many
possibilities for $\cR_{i ( \gO )}$.

So assume that $\cR \pis W$.  In particular,
$\cR \pis \cJ^W_\gO$ and $\cR$ is an element of $W$.

We claim that $\cP = \cM^{\cT}_{\eta}$.  Suppose not. 
Then $\cP$ is a fully sound proper initial segment of 
$\cM^{\cT}_{\eta}$,
and
$$A \subseteq \ran ( \gp ) \cap \gk
\subset 
H^\cR_{n+1} \left( \gr_{n+1} \left( \cP \right) 
\cup \gp \left( p_{n+1} \left( \cP  \right)  \right) \right)
\subseteq
H^\cR_\go \left( \gkbar \right)$$
where $n = n ( \cP , \gkbar )$.
But $H^\cR_\go \left( \gkbar \right)$ is an element of $W$ and has
the same cardinality as $A$.
This 
contradicts the assumption that $A$ witnesses the failure of CP$(W,\gk )$.
Therefore $\cP = \cM^{\cT}_{\eta}$.

In particular, $\eta \not= 0$, so we are done with Case~1 if $W$ satisfies
strong condensation.

Let 
$$\Cbar = \left\lbrace \crit ( E^{\cT}_\gg ) \mid
0 \leq \gg < \eta \ \land \ 
\eta \cap \cD^\cT \subseteq \gg
\right\rbrace .$$
Then
$$\cP = \cH^{\cP}_{n+1} \left( \gr_{n+1} \left( \cP \right) 
\cup \Cbar \cup p_{n+1} \left( \cP \right) \right)$$
where $n = n ( \cP , \gkbar )$.
So
$$A \subseteq \ran ( \gp ) \cap \gk
\subset 
H^\cR_{n+1} \left(\gr_{n+1} \left( \cP \right) 
\cup (\gp''\Cbar )
\cup \gp \left( p_{n+1} \left( \cP  \right)  \right) \right)$$
$$\subseteq H^\cR_\go \left( \gkbar \cup (\gp'' \Cbar ) \right)$$

Suppose for contradiction that $\Cbar$ is bounded in $\gkbar$.
Equivalently, $\gp'' \Cbar$ is bounded in $\gk$.
Since $\gk$ is least such that CP$(W , \gk )$ fails,
$$\mathrm{CP}(W , \sup ( \gp'' \Cbar ) )$$ holds. 
So there is a $B \in W$ such that
$\go_1 \cup (\cup x ) \cup (\gp '' \Cbar ) \subseteq B$ and
$$\card ( B ) \leq \card ( \go_1 \cup (\cup x ) \cup (\gp '' \Cbar ) )
\leq \card (A).$$
But then
$H^\cR_\go \left( \gkbar \cup  B \right)$
is a set in $W$ containing $A$,
and $H^\cR_\go \left( \gkbar \cup  B \right)$ has the same cardinality as $A$.
This contradicts that $A$ is a counterexample to CP$(W ,\gk )$.
Therefore $\gkbar = \sup ( \Cbar )$.

It follows from the fact that $\cT$ is thorough that $\eta$ is a limit
ordinal, $\Cbar$ is a set of indiscernibles for 
$\cJ^{\cP}_{\gkbar} = \Wbar \parallel \gkbar$,
and the last extender of $\cP$ is a measure generated by $\Cbar$.  By the
elementarity of the ultrapower map from $\cP$ into $\cR$ and its
agreement with $\gp \res \gkbar$, we conclude that clause~(b) of Theorem~5.3
holds with $C = \gp''\Cbar$.

\begin{proof}[Proof of Claim~5.3.1]
Assume for contradiction that Claim~5.3.1 is false.  So for every $i \in
S_0$, either $\Phi ( W , \cR_i )$ fails or $\cR_i \npis W$.

Let $S$ be a stationary subset of $S_0$ on which 
the choice functions $i \mapsto \eta_i$ and $i \mapsto n ( \cP_i , \gk_i )$
are constant.  Such a set $S$ exists by Fodor's lemma
\cite[Lemma 1.1]{MiSch}.  Say $\eta
= \eta_i$ and $n ( \cP_i , \gk_i ) = n$ whenever $i \in S$.

Let $j \in \lim (S) \cap S$.  Let $\psi : M \lra V_{\gO +2}$ be an
elementary embedding with $M$ countable and transitive, and everything
relevant in the range of $\psi$.  Since $\psi$ is countable, we may fix 
$i \in j \cap S$ such that 
$$\ran (\psi ) \cap X_j \subset X_i \ .$$
Let $\gp_{i,j} : N_i \lra N_j$ be the uncollapse of $(\gp_j )^{-1} (X_i)$.

We claim that $\ran ( \psi ) \cap \gk$ is unbounded in $\gk$.  The proof is
essentially the same as that of Claim~3.1.1, and we leave the details to the
reader.

Therefore 
$$\sup \left( ( \gp_{i,j} )'' \gk_i \right) = \gp_{i,j} ( \gk_i ) =
\gk_j \ .$$
Arguing as in \cite[Lemma 2.1]{MiSch}, we find a ``pull-back'' $\cP^*_j$ of $\cP_j$ by
$\gp_{i,j}$.
That is, we find an $x$-premouse $\cP_j^*$
such that $\cP^*_j$ agrees with $W_i$ below $\gk_i$ and
$$\cP_j = \ult ( \cP^*_j , \cE (\gp_{i,j} ) \res \gk_j ).$$
Moreover,
$\cP^*_j$ is $\gk_i$-sound with $n ( \cP^*_j , \gk_i ) = n$.

Since both $\cP_i$ and $\cP^*_j$ are 
$\gk_i$-sound and
iterable above $\gk_i$, one must be an initial segment of the other.
Let $\cP$ be the shorter of $\cP_i$ and $\cP^*_j$.
It will be important that $\cP \in X_j$ and that $\cP \is \cP^*_j$.

Now $\ult ( \cP , \cE ( \gp_i ) \res \gk )$ is ``bad'' in the sense that
$$\Phi \left( \cJ^W_\gO , \ult ( \cP , \cE ( \gp_i ) \res \gk ) \right)$$
fails and $\ult ( \cP , \cE ( \gp_i ) \res \gk ) \npis \cJ^W_\gO$.

Therefore,
$$N_j \models
\neg \Phi \left( 
( \gp_j )^{-1} ( \cJ^W_\gO ) ,
( \gp_j )^{-1} ( \ult ( \cP , \cE ( \gp_i ) \res \gk ) )
\right)$$
and
$$N_j \models
( \gp_j )^{-1} ( \ult ( \cP , \cE ( \gp_i ) \res \gk ) )
\npis
( \gp_j )^{-1} ( \cJ^W_\gO ) .$$
Hence,
$$N_j \models
\neg \Phi \left( W_j ,
\ult ( \cP , \cE (\gp_{i,j} \res \gk_j ))
\right)$$
and
$$N_j \models
\ult ( \cP , \cE (\gp_{i,j} \res \gk_j ))
\npis
W_j .$$
By absoluteness,
$\Phi \left(
W_j ,
\ult ( \cP , \cE (\gp_{i,j} \res \gk_j )) \right)$ fails
and
$$\ult ( \cP , \cE (\gp_{i,j} \res \gk_j ) )
\npis
W_j.$$

Now $\ult ( \cP , \cE (\gp_{i,j} \res \gk_j ))$ is iterable above $\gk_j$
since there is a sufficiently elementary embedding
from $\ult ( \cP , \cE (\gp_{i,j} \res \gk_j ))$
into an initial segment of
$$\cP_j = \ult ( \cP^*_j , \cE (\gp_{i,j} \res \gk_j )),$$
and the critical point of the embedding is at least $\gk_j$.

Moreover,
$\ult ( \cP , \cE (\gp_{i,j} \res \gk_j ))$ is $\gk_j$-sound and agrees with
$W_j$ below $\gk_j$.  In the standard way,
from the fact that we can compare
$$\ult ( \cP , \cE (\gp_{i,j} \res \gk_j ))$$ 
and $W_j$, 
we conclude that either $\Phi \left(
W_j , \ult ( \cP , \cE (\gp_{i,j} \res \gk_j )) \right)$ holds or 
$\ult ( \cP , \cE (\gp_{i,j} \res \gk_j ))
\pis W_j$.
This is a contradiction.
\end{proof}

\bigskip

\noindent
{\bf Case~2.}
{\it There is a stationary set
$$S_0 \subseteq 
\left\{ i < \ge \mid \cf (i ) > \go \ \land \ i=\ge_i \right\}$$ 
such that $i \in S_0 $ implies $(\ge_i^+ )^{W_i} = (\ge_i^+ )^W$.}

\bigskip

In this case, we conclude that (a) holds by an argument just like that given
in Case~2 of the proof of Theorem~3.2.
\end{proof}

\end{document}